\documentclass[a4paper,12pt]{article}

\usepackage{amsthm,amsmath,amssymb,latexsym}
\numberwithin{equation}{section}
 
\newcommand{\QQ}{{\mathbb Q}}

\newcommand{\Sing}{\mathop{\rm Sing}\nolimits}

\begin{document}

\begin{center}
\large{\textbf{Uniserial Noetherian Centrally Essential Rings}}
\end{center}

\hfill {\sf V.T. Markov}

\hfill Lomonosov Moscow State University

\hfill e-mail: vtmarkov@yandex.ru

\hfill {\sf A.A. Tuganbaev}

\hfill National Research University "MPEI"

\hfill Lomonosov Moscow State University

\hfill e-mail: tuganbaev@gmail.com

{\bf Abstract.} It is proved that a ring $R$ is a right uniserial, right Noetherian centrally essential ring if and only if $R$ is a commutative discrete valuation domain or a left and right Artinian, left and right uniserial ring. It is also proved that there exist non-commutative uniserial Artinian centrally essential rings.\footnote{V.T.Markov is supported by the Russian Foundation for Basic Research, project 17-01-00895-A. A.A.Tuganbaev is supported by Russian Scientific Foundation, project 16-11-10013.}


\section{Introduction}
All rings considered are associative and contain a non-zero identity element. Writing expressions of the form ``$R$ is an Artinian ring", we mean that the both modules $_RR$ and $R_R$ are Artinian.

A ring $R$ with center $C$ is said to be \emph{centrally essential} if  for any non-zero element $r\in R$, there exist two non-zero central elements $c,d\in C$ with $rc=d$\footnote{This is equivalent to the property that the module
$R_C$ is an essential extension of the module $C_C$.}
  
Centrally essential rings are also studied in \cite{MT18}, \cite{MT18b}, \cite{MT18c}, \cite{MT18d}, \cite{MT19a}, \cite{MT19b}, \cite{MT19c}. 

Since any centrally essential ring, which is semiprime or right non-singular ring, is com\-mu\-ta\-ti\-ve (see \cite[Proposition 3.3]{MT18} and \cite[Proposition 2.8]{MT19a}), it is interesting to study only centrally essential rings which are not semiprime or right non-singular.

We note that the exterior algebra $\Lambda(V)$ of a finite-dimensional vector space $V$ over a field $F$ of characteristic $0$ or $p\neq 2$ is centrally essential if and only if
$\dim V$ is an odd positive integer; in particular, if $F$ is a finite field of odd characteristic and $\dim V$ is an odd positive integer exceeding $1$, then $\Lambda(V)$ is a centrally essential noncommutative finite ring; see \cite{MT18c}. Consequently, centrally essential PI rings are not necessarily commutative.

We also note that if $Z_2$ is the field of order 2 and $Q_8]$ is the quaternion group of order 8, then the group ring $Z_2[Q_8]$ is a local finite centrally essential ring which is not commutative, \cite{MT18}.

It is well known that a commutative domain $R$ is a uniserial Noetherian ring (resp., a uniserial ring) if and only if $R$ is a discrete valuation domain (resp., a valuation domain).

It follows from the above that uniserial Noetherian centrally essential prime rings (resp., uniserial centrally essential prime rings) coincide with commutative discrete valuation domains (resp., commutative valuation domains).

The following theorem is the main result of this paper.

\textbf{Theorem 1.}\label{main}
A centrally essential ring $R$ is a right uniserial, right Noetherian ring if and only if $R$ is a commutative discrete valuation domain or a (not necessarily commutative) uniserial Artinian ring.

\textbf{Remark.} In connection to Theorem 1, we remark that there exist right uniserial right Noetherian rings which are neither prime rings nor right Artinian rings; e.g., see \cite[Example 9.10(3)]{Tug98}.

The proof of Theorem 1 is contained in the next section and is based on several assertions, some of which are of independent interest.

We give some notation and definitions. A module $M$ is said to be \emph{uniserial} if the set of all submodules of the module $M$ is linearly ordered with respect to inclusion. A module $M$ is said to be \textit{uniform} if any two non-zero submodules of $M$ have the non-zero intersection. A module $M$ is called an \textit{essential extension} of some its submodule $X$ if $X\cap Y\ne 0$ for any non-zero submodule $Y$ of $M$. For a right $R$-module $M$, we denote by $\Sing M$ the \textit{singular submodule} of $M$, i.e., $\Sing M$ is the submodule of $M$ consisting of all elements whose right annihilators are essential right ideals of $R$.

For a ring $R$, we denote by $J(R)$, $\Sing R_R$, $R^*$ and $C(R)$ the Jacobson radical, the right singular ideal, the group of invertible elements and the center of the ring $R$, respectively. The \emph{left annihilator} of an arbitrary subset $S$ of is $\ell .Ann_R(S)=\{r\in R\,|\, rS=0\}$. The right annihilator $r.Ann_R(S)$ of the set $S$ is similarly defined.
We denote by $[a,b]=ab-ba$ the \emph{commutator} of the elements $a,b$ of an arbitrary ring; we also use the following well known obvious properties of commutators: $[a,b]=-[b,a]$, $[ab,c]=a[b,c]+[a,c]b$ for any three elements $a,b,c$ of an arbitrary ring.

Other ring-theoretical notions and notations can be found in \cite{Herstein, Lambek, Tug98}.

\section{The proof of Theorem 1}\label{section2}

\textbf{Lemma 2.1.}\label{zd-ideal}
If the set $B$ of all left zero-divisors of the ring $A$ is an ideal, then
$B$ is a completely prime ideal.

\textbf{Proof.} Let $a,b\in A$ and $ab\in B$. Then there exists an element $x\in A\setminus\{0\}$ with $abx=0$. If $bx=0$, then $b\in B$. Otherwise, it follows from the relation $a(bx)=0 $ that $a\in B$.~\hfill$\square$

\par\noindent\textbf{Lemma 2.2.} \label{zero-div} Let $A$ be a centrally essential ring. Then:

\textbf{2.2.1.} in $A$, all one-sided zero-divisors are two-sided zero-divisors;

\textbf{2.2.2.} the ring $A$ is left uniform if and only if $A$ is right uniform;

\textbf{2.2.3.} if the ring $A$ is right uniform and $B=\Sing A_A$, then $B$ is
the set of all (left or right) zero-divisors of the ring $A$ and $B$ is a completely prime ideal of the ring $A$;

\textbf{2.2.4.} if the ring $A$ has a proper ideal $B$ containing all left zero-divisors of the ring $A$, then the factor ring $A/B$ is commutative;

\textbf{2.2.5.} if an ideal $B$ of the ring $A$  contains all central zero-divisors of the ring $A$, then $\ell .Ann_A(B)\subseteq{}C(A)$.

\textbf{Proof.}

\textbf{2.2.1.} Let $a_1,a_2$ be two non-zero elements of the centrally essential ring $A$ with $a_1a_2=0$. There exist non-zero central elements $x_1,x_2,y_1,y_2$ of $A$ such that $a_1 x_1=y_1$, $a_2 x_2=y_2$. Then $y_2a_1 =a_1y_2 =a_1a_2x_2=0$. Therefore, the left zero-divisor $a_1$ is a right zero-divisor. Similarly, the right zero-divisor
$a_2$ is a left zero-divisor.
 
\textbf{2.2.2.} Let's assume that the ring $A$ is right uniform and $a_1,a_2$ are two non-zero elements of the ring $A$. There exist non-zero central elements
$x_1,x_2,y_1,y_2$ of $A$ such that $a_1 x_1=y_1$, $a_2 x_2=y_2$. Then $$Aa_1\cap Aa_2\supseteq Ax_1a_1\cap Ax_2a_2=a_1x_1A\cap a_2x_2A=y_1A\cap
 y_2A\neq 0.$$
 
\textbf{2.2.3.} By the definition of the right singular ideal, all its elements are left
zero-divisors. Conversely, let $a$ be a left or right zero-divisor of the ring $A$. Then $r.Ann(a)\neq 0$ by the first assertion of the lemma; in a right uniform ring, this means that $r.Ann(a)$ is an essential right ideal, i.e.,  $a\in B$. Now we use Lemma 2.1.\\

\textbf{2.2.4.} Let $a,b\in A\setminus B$. There exist two non-zero central elements $x,y$ of $A$ with $b x=y$. Then $[a,b]x=[a,bx]=0$, i.e., $[a,b]$ is a left zero-divisor. Consequently, $[a,b]\in B$.

\textbf{2.2.5.} Let $r\in \ell .Ann_A(B)$. There exist two non-zero central elements $x,y$ of $A$ with $rx=y$. It is clear that $x\not\in B$; therefore, $x$ is not a zero-divisor. Therefore, for every element $a\in A$, it follows from the relations $0=[a,y]=[a,rx]=[a,r]x$ that $[a,r]=0$.~\hfill$\square$

\textbf{Corollary 2.3.} In a right uniform centrally essential ring, the left singular ideal coincides with the right singular ideal.

For convenience, we give brief proofs of the following two well known assertions.

\textbf{Lemma 2.4.}\label{div-fg-mod} Let $A$ be a commutative domain which has a non-zero finitely  generated divisible torsion-free $A$-module $M$. Then $A$ is a field.

\textbf{Proof.} Let's assume the contrary. Then $A$ has a non-zero maximal ideal $\mathfrak m$ and $M$ naturally turns into a non-zero finitely  generated module over the local ring $R_{\mathfrak m}$ with radical $J={\mathfrak m}R_{\mathfrak m}$. Since the module $M$ is divisible, we have that $MJ\supseteq M{\mathfrak m}=M$ and $M=0$ by the Nakayama lemma. This is a contradiction.~\hfill$\square$

\textbf{Lemma 2.5.}\label{N=aN}
Let $A$ be a right uniserial ring and $P$ a completely prime ideal of $A$. Then $N=aN$
for every $a\in A\setminus P$.

\textbf{Proof.} Let $a\in A\setminus P$. Since $aA\not\subseteq P$, we have $P\subseteq aA$. Therefore, for every $x\in P$, there exists an element $b\in A$ with $x=ab$. Since $P$ is a completely prime ideal, $b\in P$ and $x\in aP$.~\hfill$\square$

\textbf{Proposition 2.6.}\label{right-artinian} Let $A$ be a right uniserial, right Noetherian, centrally essential ring. Then $A$ is either a commutative domain or a right uniserial, right Artinian ring.

\textbf{Proof.} We set $N=\Sing A_A$. The ideal $N$ is nilpotent; e.g.,  see \cite[9.2]{Tug98}. It follows from Lemma 2.2(3,4) that the ideal $N$ is completely prime and  contains all zero-divisors of the ring $R$ and the ring $A/N$ is a commutative domain. Therefore,  the proposition is true for $N=0$. Now let $N\neq 0$. We denote by $n$ the nilpotence index of the ideal $N$. Then $0\neq N^{n-1}\subseteq{}\ell .Ann_A(N)$. It follows from Lemma 2.2(5) that $N^{n-1}\subseteq C(A)$. Next, for every $a\in A\setminus N$, we have $N=aN$ by Lemma 2.5, whence $N^{n-1}=aN^{n-1}=N^{n-1}a$. Consequently, $N^{n-1}$ is a divisible right $(A/N)$-module and $N^{n-1}$ is a torsion-free $(A/N)$-module since all zero-divisors of the ring $A$ are contained in $N$. By Lemma 2.4, the ring $A/N$ is a field and each of the cyclic {(A/N)}-modules $(N^{k-1}/N^k)$ for $k=1,\ldots,n$ is a simple module. Consequently, the ring $A$ is right Artinian.~\hfill$\square$

\textbf{Lemma 2.7.}\label{local-pri}
Let $A$ be a local ring and let $J(A)=\pi A$ for some element $\pi\in A$ of nilpotence index $n$ (maybe, $n=\infty$). For any two integers $k,\ell $ and each $a,b\in A$ such that $k,\ell \geq 0$, $k+\ell <n$, $a\in\pi^kA\setminus\pi^{k+1}A$ and $b\in\pi^lA\setminus\pi^{\ell +1}A$, we have $ab\in\pi^{k+\ell }A\setminus\pi^{k+\ell +1}A$.

\textbf{Proof.} It follows from the inclusion $A\pi\subseteq \pi{}A$ that $ab\in \pi^{k+\ell }A$. If $\pi^m\in \pi^{m+1}A$ for some $m\geq 0$, then it is clear that $\pi^m(1-\pi t)=0$ for some $t\in A$ and $\pi^m=0$, since $1-\pi t\in A^*$. We set $a=\pi^{k}r$ and $b=\pi^{\ell }s$ for some $r,s\in A\setminus J(A)$. Then $r,s\in A^*$, since the ring $A$ is local and $r\pi^\ell \in\pi^lA\setminus\pi^{\ell +1}A$. Consequently, $r\pi=\pi{}r'$ for some $r'\in A^*$ and $ab=\pi^{k+\ell }r's$. It remains to remark that $ab\not\in\pi^{k+\ell +1}A$, since $r's\in A^*$ and
$\pi^{k+\ell }\neq 0$.~\hfill$\square$

\textbf{Proposition 2.8.}\label{left} A right uniserial, right Artinian, centrally essential ring is a left uniserial, left Artinian ring.

\textbf{Proof.}
Let $A$ be a right uniserial, right Artinian, centrally essential ring, $N=J(A)$, and let $n$ be the nilpotence index of the ideal $N$. If $n=1$, then the ring $A$ is commutative by Lemma 2.2(3,4); it is nothing to prove in this case. Any right uniserial ring is a local ring; therefore, every element of $A\setminus N$ is invertible. Let $n>1$, i.e., $N\neq 0$. Since a right Noetherian (e.g., a right
Artinian) right uniserial ring is a principal right ideal ring, $N=\pi A$ for some element $\pi\in N$. There exist two elements $x,y\in C(A)$ with $\pi x=y\neq 0$.
Let $x\in N^k \setminus N^{k+1}$ for some $k$, $0\leq k<n$. Then $y\in N^{k+1}$, whence $k+1<n$. If $[a,\pi]\not \in N^2$, then it follows from Lemma 2.7 that $[a,\pi]x\not\in N^{k+2}$; consequently, $[a,\pi]x\neq 0$. However, $[a,\pi]x=[a,\pi x]=[a,y]=0$.
This is a contradiction; therefore, $[a,\pi]\in N^2$ for every $a\in A$.
Consequently, $N=A\pi+N^2$, whence $N/A\pi=N(N/a\pi)=\ldots=N^n(N/A\pi)=0$, i.e., $N=A\pi$.
It follows from the left-side analogue of Lemma 2.7 that every left ideal of the ring $A$ coincides with one of the ideals $A,N,N^2, \ldots, N^{n-1},\{0\}$, i.e., $A$ is a left uniserial, left Artinian ring.~\hfill$\square$

\textbf{Proposition 2.9.}\label{exists}
Let $F$ be a field and let $D_1,D_2\colon F\rightarrow F$ be two derivations of the field $F$ with incomparable kernels (for example, we can take the field of rational functions $\QQ(x,y)$ in two independent variables as $F$ and set $D_1=\partial/\partial x$, $D_2=\partial/\partial y$).

Then for every positive integer $n\geq 2$, there exists a non-commutative uniserial, Artinian, centrally essential ring $A$ such that $A/J(A)\cong F$ and the nilpotence index of the ideal $J(A)$ is equal to $n$.

\textbf{Proof.} We use a construction which is similar to the one described in \cite{Jelis}. Let $N=2n-1$, $R=M_N(F)$ be the matrix ring of order $N$ over the field $F$, $e_{i,j}$ denote the matrix unit for any $i,j\in\{1,\ldots,N\}$, and let $f\colon F\rightarrow R$ be the mapping defined by the rule
$$
f(\alpha)=\alpha E+D_1(\alpha)e_{1,N-1}+D_2(\alpha)e_{N-1,N}
$$
for every $\alpha\in F$, where $E$ is the identity matrix. Let $A$ be the subring of the ring $R$ generated by the set $f(F)$ and the matrix
$\pi=\sum_{i=1}^{n-1}e_{2i-1,2i+1}$. It is directly verified that 
$\pi^n=0$, $\pi^{n-1}=e_{1,N}$, $f(\alpha)\pi=\pi f(\alpha)=\alpha \pi$ and 
$$[f(\alpha),f(\beta)]=(D_1(\alpha)D_2(\beta)-D_1(\beta)D_2(\alpha))\pi^{n-1}$$ 
for any $\alpha,\beta\in F$. It follows from these relations that $\pi A=A\pi=J(A)$, $J(A)^k=\pi^kA=A\pi^k$ for all $k=1,\ldots,n-1$ and $\pi A\subseteq C(A)$. It is clear that $A$ is a uniserial Artinian ring. If $a\in A\setminus \{0\}$ and $a\in \pi A$, then $a\in C(A)$; otherwise, $a\pi^{n-1}\in C(A)\setminus \{0\}$ and $X^{n-1}\in C(A)$. Consequently, the ring $A$ is centrally essential.

Finally, if $\alpha\in \ker D_2\setminus \ker D_1$ and $\beta\in \ker D_1\setminus \ker D_2$, then $$[f(\alpha),f(\beta)]=D_1(\alpha)D_2(\beta)X^{n-1}\neq 0,$$ i.e., the ring $A$ is not commutative.~\hfill$\square$

\textbf{2.10. Completion of the proof of Theorem 1.} The first assertion follows from Propositions 2.6 and 2.8. The second assertion follows from Proposition 2.9.

\textbf{2.11. Remark.} It is known that if $R$ is a right noetherian right uniserial ring and $\bigcap_{n\ge0}J(R)^n\ne 0$, then $R$ is not left noetherian and $R$ is not left Ore; see for instance Proposition 3.7 in [1]. As a corollary, $\bigcap_{n\ge0}J(R)^n=0$ for every noetherian right uniserial ring. Also,  $\bigcap_{n\ge0}J(R)^n=0$ for every uniserial right noetherian ring. These results are related to Theorem 1, and to several results in Section 2.

The authors are grateful to the reviewer for useful comments.

\end{document}